\documentstyle[12pt,amscd, amssymb]{amsart}

\newcommand{\cI}{{\cal I}}
\newcommand{\cO}{{\cal O}}

\newcommand{\cx}{{\cal X}}
\newcommand{\cH}{{\cal H}}
\newcommand{\cZ}{{\cal Z}}

\newcommand{\ol}{\overline}
\newcommand{\pn}[1]{{\Bbb P}^{#1}}

\newcommand{\D}{\Delta}
\newcommand{\s}{\sigma}
\newcommand{\e}{\epsilon}
\newcommand{\ra}{\rangle}
\newcommand{\la}{\langle}

\newcommand{\pr}[2]{{\partial}_{#1}^{#2}}
\newcommand{\fm}{{\frak m}}

\newtheorem{theorem}{{\rm T\sc heorem}}[section]
 \newtheorem{lemma}[theorem]{{\rm L\sc emma}}
 \newtheorem{corollary}[theorem]{{\rm C\sc orollary}}
\newtheorem{proposition}[theorem]{{\rm P\sc roposition}}

\newtheorem{definition}[theorem]{{\rm D\sc efinition}} 
\newtheorem{re}[theorem]{{\rm R\sc emark}}

\theoremstyle{remark}
\newtheorem{rem}[theorem]{{\rm R\sc emark}}
\newtheorem{ex}[theorem]{\rm E\sc xample}

\newenvironment{rem*}{\begin{re}\em}{\end{re}}
\newenvironment{paragraph*}[1]{\indent\begin{paragraph}{\bf #1.}\em\
}{\end{pargrph}}

\begin{document}


\title{ Generation of $k$-jets on Toric Varieties}

\author{Sandra Di Rocco}
\address{\small KTH, 
Department of Mathematics \\          
 S-100 44 Stockholm}\email{ sandra@@math.kth.se}
\curraddr{Max-Planck-Institute Fur Mathematik\\
Gottfried-Claren-Str. 26
D-53225 Bonn}
\subjclass{ Primary 14M25, 14J60, 14C20;
Secondary 14C25, 14E25}
\keywords{Toric varieties,
convex polyhedra, $k$-jet ampleness}
\thanks{Partially supported by NFR B/276900}

\maketitle

\begin{abstract} The notion of a $k$-convex $\D$-support function  for a toric variety $X(\D)$ is introduced. A criterion for a line bundle $L$ to generate $k$-jets on $X$ is given in terms of the $k$-convexity of the $\D$-support function $\psi_L$. Equivalently $L$ is proved to be $k$-jet ample if and only if the restriction at each invariant curve has degree at least $k$.
\end{abstract}

\section*{Introduction}
The notion of $k$-jet ampleness has been introduced by Demailly to describe line bundles $L$  whose global sections can have arbitrarily prescribed $k$-jets at every single point $x\in X$, see \cite{Dem}. Beltrametti and Sommese generalized it by considering $k$-jets supported on a finite number of points.
\vskip6pt
{\em A line bundle $L$ is said to be $k$-jet ample on $X$ if for any collection of $r$ points, $(x_1,\cdots,x_r)$, and any $r$-ple of positive integers $(k_1,\cdots,k_r)$, with $\sum k_i=k+1$ the natural map
$$H^0(X,L)\to H^0(L\otimes\cO_X/\fm_{x_1}^{k_1}\otimes \cdots \otimes \fm_{x_r}^{k_r})$$
is surjective}.
\vskip6pt
 Notice that $0$-jet ampleness is equivalent to being spanned by global sections and $1$-jet ampleness is equivalent to being very ample.

During the last years many results on $k$-jet ample line bundles on surfaces have been established, \cite{BeSok, EiLa, BaSz, BaDRSz}. Up to our knowledge the problem is still quite open for higher dimensional varieties, besides few cases like $\pn{n}$ and Fano varieties \cite{BeSok, BeDRSo}.

In \cite{Cox1} D. Cox has introduced ``homogeneous coordinates" on a toric variety $X(\D)$. For the points invariant under the torus action the situation looks  similar to the projective space case. Using this system of local coordinates we give a description of the fibers of the $k$-jet bundle $J_k(L)$ on fixed points, see section \ref{kjet}. 

According to Oda and Demazure a line bundle $L$ on $X(\D)$ is generated by global sections (respectively very ample) if the $\D$-support function $\psi_L$ is {\em convex (respectively strictly convex)}. This suggests to use a ``higher convexity" property for $\psi_L$ in the cases $k\geq 2$.\\
In section \ref{kconvex} we introduce the notion of a {\em $k$-convex $\D$-support function}, which for $k=0,1$ agrees with being convex or strongly convex.\\
The $\D$-support function $\psi_L$ is $k$-convex if the polyhedra $P_L$, associated to $L$,  has edges of length at least $k$. This translates to the property that the intersection of $L$ with the invariant curves, associated to every edge, is $\geq k$, which is a generalization of the  toric Nakai criterion for ample line bundle.

A key step in the proof is the reduction to the case where the considered points are invariant under the torus action. We are grateful to T. Ekedahl for suggesting to use the Borel's fixed point theorem, and for pointing out the sufficiency of the reduction argument.

Our result states that in order to check the $k$-jet ampleness of a line bundle $L$ it is enough to have a bound on the intersection $L\cdot C$, for all the invariant curves $C$. This can be applied to the study of ``local positivity". In section \ref{applications} we report a series of results on blow-ups and higher adjoint bundles, which in the toric case can be shown by means of a direct checking on intersections. We also state an equivalent criterion for $k$-jet ampleness in terms of a bound of the Seshadri constant $\e(L,x)$. This can be thought as a toric version of the Seshadri criterion for ample line bundles, generalized to $k$-jet ampleness.

In this paper we prove:\\
{\em Let $L$ be a line bundle on a non singular toric variety $X(\D)$, then the following statements are equivalent:
\begin{itemize}
\item $L$ is $k$-jet ample;
\item $L\cdot C\geq k$, for any $T$-invariant curve $C$, [Proposition \ref{inter}];
\item $\psi_L$ is $k$-convex, [Theorem \ref{main}];
\item  the Seshadri constant $\e(L,x)\geq k$ for each $x\in X$, [Proposition \ref{sesh}].
\end{itemize}}
\vskip6pt

It is a pleasure to acknowledge valuable discussions with T. Ekedahl, M. Boij, D. Laksov and D. Cox. This research was initiated during the author's stay at Mittag-Leffler Institute and KTH and reached a final stage at the Max Planck Institute. To all those institutions we owe thanks for their support.

\section*{ Notation}
We will use standard notation in Algebraic Geometry. The groundfield will always be the field of complex numbers.\\
By abuse of terminology the words line bundle and Cartier divisor will be used with no distinction, as well as the multiplicative and additive structure.\\
For basic notions on toric varieties we refer to \cite{Fu, Oda, Ew} and for a nice survey on the resent progress on toric geometry we refer to \cite{Cox2}.\\
When not stated $X$ will always denote a smooth $n$-dimensional toric variety and  $L$ a line bundle on it.\\
\section{Toric Varieties}
Let $N$ be an $n$-dimensional lattice and $\D=\cup \s_i$ be a complete and regular fan, meaning: 
\begin{itemize}
\item $supp(\D)=N_{\Bbb R}=N\otimes {\Bbb R}$ and
\item for every  $r$-dimensional cone $\s\in\D$, there exists a $\Bbb Z$-basis of $N_{\Bbb R}$, $\{\rho_1,\cdots,\rho_n\}$, such that the subset $\{\rho_1,\cdots,\rho_r\}$ spans $\s$.
\end{itemize}
 We will denote by $X=X(\D)$ the associated  non singular $n$-dimensional toric variety and by $\D(t)$ the set of $t$-dimensional cones in $\D$.\\
Let $M=Hom_{\Bbb Z}(N,{\Bbb Z})$ be the dual lattice so that $X$ is obtained by gluing together the affine toric varieties $X_{\s}=Spec({\Bbb C}[\check{\s}\cap M])$, where $\check{\s}=\{v\in M_{\Bbb R}:\la v,\s\ra\geq 0\}$, and $\s\in\D$.\\
Each $m\in M$ can be viewed as a rational function $\chi^m:T=N\otimes{\Bbb C}^*\to \Bbb C^*.$\\
There is a $1-1$ correspondence between $r$-dimensional cones $\tau\in\D(r)$ and $T$-invariant codimension $r$ subvarieties of $X$, which will be denoted by $V(\tau)$.\\
Let $D_i=V(\rho_i)$ be the $T$-invariant divisors corresponding to the one dimensional cones $\rho_i\in\D(1)$.  The set $\{D_i\}_{\rho_i\in\D(1)}$ form a set of generators for the Picard group of $X$ and thus every line bundle $L$ can be written in terms of those principal divisors: 
$$L=\sum_{\rho_i\in\D(1)} a_i D_i$$
We will denote by $P_L$ the associated convex polyhedra:
$$P_L=\{m\in M_{\Bbb R}:\la m, \rho_i\ra\geq -a_i\}$$
This gives a nice way of expressing the global sections of $L$:
$$H^0(X,L)=\bigoplus_{m\in P_L\cap M}{\Bbb C}\cx^m$$
Recently David Cox \cite{Cox1}, has introduced the notion of homogeneous coordinates on a toric variety.
There is a $1-1$ correspondence between $T$-invariant principal divisors $D_i$ and linear monomials $\cx_i$ on $X$.
The polynomial ring is then defined as:
$$S={\Bbb C}[\cx_i:\rho_i\in\D(1)]$$
and the grading is given by the group of divisors modulo rational equivalence, $Pic(X)$, i.e. two rationally equivalent divisors $D$ and $E$ are associated to monomials $\cx_D$ and $\cx_E$ of the same degree.\\
Considering the exact sequence:
$$0\to M\to \oplus_{\rho_i\in\D(1)}{\Bbb Z}\cdot D_i\to Pic(X)\to 0$$
we associate to each $m\in M$ a divisor $\sum \la m,\rho_i\ra D_i=div(\chi^m)$.\\
 The global sections of $L$ are generated by the monomials of the form $$(\Pi_i \cx_i^{\la m,\rho_i\ra+a_i})_{m\in P_D\cap M}$$
The notion of $\D$-support function will be use constantly throughout this paper:
 \begin{definition}\cite[2.1]{Oda} A real valued function 
$f:\cup_i\sigma_i\to\bf R$ is a $\D$-linear support function if it is
$\Bbb Z$-valued on $N\cap(\cup_i\sigma_i)$ and it is linear on each
$\sigma_i$.
\end{definition} This means that for each $\sigma$ there exists $m_\sigma\in M$
such that
$f(n)=\langle m_\sigma,n\rangle$ for $n\in\sigma$ 
and  $\langle m_\sigma,n\rangle=\langle m_\tau,n\rangle$ when $\tau$ is a face of $\s$.
To each divisor $L$ we associate a $\D$-support function $\psi_L$ defined by:
$$\psi_L(\rho_i):=-a_i$$

\section{$k$-jet Bundles}\label{kjet}
Let $\Delta$ be the diagonal in $X\times X$ and $p:X\times X\to X$ the projection onto the first factor. The $k$-th jet-bundle associated to $L$ is the vector bundle  associated to the sheaf:
$$p^*L/p^*L\otimes \cI_{\Delta}^{k+1}$$ 
Where $\cI_{\Delta}$ is the ideal sheaf of $\Delta$. It is a vector bundle of rank ${k+n}\choose n$ whose fiber is
$$J_k(L)_x=L_x\otimes\cO_X/m_x^{k+1}$$
For details on jet bundles we refer the reader to \cite[Ch.I]{KuSp}.
There are natural maps (defined on the sheaf level):
$$i_k:L\to J_k(L)$$
sending the germ of a section $s$ at a point $x\in X$ to its $k$-th jet.
More specifically for $s\in H^0(X,L)$ $i_k(s(x))\in \bigoplus_1^{{k+n}\choose n}\Bbb C$ is the ${k+n}\choose n$-ple determined by the coefficients of the terms of degree up to $k$, in the Taylor expansion of $s$ around $x$.\\
So if $(x_1,\cdots,x_n)$ are local coordinates around $x_0=(0,0,\cdots,0)$ and $s=\sum c_{i_1,\cdots,i_r}\prod x_i^{i_j}$ then 
$$i_k(s(x_0))=(\cdots,\frac{\partial s}{\pr{x_1}{t_1}\cdots\pr{x_r}{t_r}},\cdots)|_{x=x_0}=(\cdots,(\text{{\Tiny constant}})\cdot c_{t_1,\cdots,t_r},\cdots)  $$
where $t_1+\cdots+t_r\leq k$. For example $i_1(s(x))$ consists in the constant and linear term.
The following definition formalizes the property for a linear series $|L|$ on $X$ to generate $k$-jets on one or more points of $X$. When more points are considered $L$ is said to generate ``simultaneous jets" at those points. 
  
\begin{definition}
Let $\cZ=\{x_1,\cdots,x_r\}$ be a finite collection of distinct points on $X$. $L$ is said to be $k$-jet ample on $\cZ$ (or equivalently the series $|L|$ is said to generate all $k$-jets on $\cZ$) if for any $r$-ple of positive integers $(k_1,\cdots,k_r)$, such that $\sum_1^r k_i=k+1$ the map:
$$H^0(X,L)\to H^0(L\otimes \cO_X/{\frak m}_{x_1}^{k_1}\otimes\cdots\otimes {\frak m}_{x_r}^{k_r})$$ is surjective, where ${\frak m}_{x}^k$ is $k$-th tensor power of the maximal ideal sheaf ${\frak m}_x$. 
$L$ is $k$-jet ample on $X$ if it is $k$-jet ample on each such $\cZ$ in $X$.
\end{definition}
Clearly from the definition:
\begin{itemize}
\item We can rewrite the map above as:
$$\psi_{\cZ}^{k_1,\cdots,k_r}:H^0(X,L)\to \bigoplus_1^r (J_{k_i-1}(L))_{x_i}$$
defined by $\psi_{\cZ}^{k_1,\cdots,k_r}(s)=(i_{k_1-1}(s(x_1)),\cdots,i_{k_r-1}(s(x_1)))$.
We say then that $L$ is $k$-jet ample on $X$ if the map $\psi_{\cZ}^{k_1,\cdots,k_r}$ is surjective for any $\cZ$ and any $(k_1,\cdots,k_r)\in\Bbb Z^r_+$, such that $\sum k_i=k+1$.
\item  If $L$ is $0$-jet ample then $L$ is generated by its global sections;
\item if $L$ is $1$-jet ample then using the sections in $H^0(X,L)$ we can define an embedding $i:X\to \pn{N}$ and thus $L$ is very ample. 
\end{itemize}
Using the homogeneous coordinates introduced in the previous section the $k$-jets at the $T$-invariant points $x(\s)=V(\s)$ can be better described in terms of the polyhedra associated to $L$.\\
Let $\s=\la \rho_1,\cdots,\rho_n\ra$, where $\s\in\D(n)$ and $\rho_i\in\D(1)$ are the one dimensional cones generating $\s$.
The point $x(\s)$, lies on the intersection of the divisors $D_i$, $i=1,\cdots,n$: $x(\s)\in\cap_1^n(\cx_i=0)$. Then the maximal ideal is generated by the linear monomials in $\cx_i$:
$$\fm_{x(\s)}=\la \cx_1,\cdots,\cx_n\ra$$
and thus
$$\fm_{x(\s)}^{k+1}=\la \Pi_{\rho_i\subset\s}\cx_i^{t_i}|t_1+\cdots+t_n=k+1\ra$$
i.e.  the generators are the monomials of ``degree$=k+1$" in the variables $\cx_1,\cdots,\cx_n$ ( here by degree we mean  the sum of the powers of the variables, i.e. the usual one). \\
Each ${\cx}^m$, generator of $H^0(X,L)$, can be written in the local coordinates $(\cx_1,\cdots,\cx_n)$ as follows. Fix $\{\rho_1,\cdots,\rho_n\}$ as basis of $N$ and let $\{m_1,\cdots,m_n\}$ the dual basis. In this coordinate system $m=\sum \la m,\rho_i\ra m_i$ and the germ of ${\cx}^m$ at $x(\s)$ is:
$${\cx}^m|_{x(\s)}=\prod_{i=1}^n \cx_i^{\la m,\rho_i\ra+a_i}$$
Taking its $k$-th jet means ``killing" all the monomials of degree $\geq k+1$ in the variables $\cx_1,\cdots\cx_n$:
$$ i_k({\cx}^m(x(\s)))=(\cdots, \frac{\partial{\cx}^m}{\pr{x_1}{t_1}\cdots\pr{x_r}{t_r}},\cdots)|_{x=x(\s)}$$
\ex\label{delpezzo} Let $N={\Bbb Z}^2$ and $\D$ be the $2$-dimensional fan composed by the following $6$ cones, and their edges:\\
$\s_1=\la (0,1),(1,1)\ra,\,\s_2=\la (1,1),(1,0)\ra,\,\s_3=\la (1,0),(0,-1)\ra\\
\s_4=\la (0,-1),(-1,-1)\ra,\,\s_5=\la (-1,-1),(-1,0)\ra,\,\s_6=\la (-1,0),(0,-1)\ra\\$
$X(\D)$ is the equivariant blow up of $\pn{2}$ in the $3$ fixed points, i.e. a Del Pezzo surface of degree $6$.\\
Let $L=D_1+D_2+D_3+D_4+D_5+D_6=-K_{X(\D)}$, where the $D_i's$ are associated to the edges in the order given above. Let $\s=\la (0,1),(1,1)\ra$ and let $\{m_1,m_2\}$ be the basis dual to $\{(0,1),(1,1)\}$.\\
In this basis $P_L$ is the convex hull of the points $$\{(0,1),(1,1),(1,0),(-1,0),(-1,-1),(-1,0)\}$$ 
and thus the generators of $H^0(X,L)$ are  $$\{1,\cx_1,\cx_2,\cx_1\cx_2,\cx_1^2\cx_2,\cx_1\cx_2^2,\cx_1^2\cx_2^2\}$$
Moreover ${\frak m}_{x(\s)}^2=\la \cx_1\cx_2,\cx_1^2,\cx_2^2\ra$ and then
$$J_1(L)_{x(\s)}={\Bbb C}\oplus {\Bbb C}\cx_1\oplus \cx_2$$

\ex\label{pn1} Let $X=\pn{n}$, then $\D$ is the fan composed by $(n+1)$ $n$-dimensional cones spanned by the $(n+1)$ edges
\begin{itemize}
\item $\rho_i=(0\cdots,0,\underbrace{1}_{i-th},0,\cdots,0)$ for $i=1,\cdots n$
\item $\rho_{n+1}=(-1,\cdots,-1)=-rho_1-\cdots -\rho_n$
\end{itemize}
 Let $D_1,\cdots,D_{n+1}$ be the associated $T$-invariant principal divisors and let $L=D_1+\cdots+D_k=\cO_{\pn{n}}(k)$.\\
recall that the Picard group is generated by one principal divisor $D_i$ and that $D_i\equiv D_j$ for $i\neq j$. So we can think of $L$ as
$$L=t_1D_1+\cdots+t_nD_n;\,\,t_1+\cdots+t_n=k$$
Let $\s=\la\rho_1,\cdots\rho_n\ra$, and fix the basis $\{\rho_1,\cdots,\rho_n\}$ with dual $\{m_1,\cdots,m_n\}$. In this basis the polyhedra $P_L$ is the convex hall of the $(n+1)$ points
$$\{ (-1,\cdots,-1),(k,-1,\cdots,-1),\cdots,(-1,\cdots,-1,k,-1,\cdots,-1),\cdots,(-1,\cdots,k)\}$$
Then for any decomposition $t_1,\cdots,t_n$ of positive integers such that $\sum_1^n(t_1+1)=k$ the lattice point $m=\sum _1^n t_1m_1\in P_L$ and any lattice point $m\in P_L$ can be written in this form. The situation stays the same if we consider another $\s\in\D$. It follows that 
$$J_k(L)_{x(\s)}=H^0(X,L)=\bigoplus_{t_1+\cdots+t_r=k}{\Bbb C}\prod_1^r\cx_1^{t_i}$$
In particular if $L=\cO(1)$ then $J_1(L)$ is trivial. This is in fact a characterization of the projective space, cf.\cite{So}.
\section{ $k$-convex functions}\label{kconvex}
In order to study positivity properties of line bundles Demazure and Oda introduced the definition of convex and strictly convex $\D$-support function.
\begin{theorem}\cite[Th. 2.13]{Oda} A line bundle $L$ on $X$ is globally generated ( i.e. $0$-jet ample) if and only if $\psi_L$ is convex and it is very ample (i.e. $1$-jet ample) if and only if $\psi_L$ is strictly convex.
\end{theorem}
A natural way of generalizing such a criterion to higher jets is to introduce a definition of ``higher convexity".
 
 \begin{definition}  Let $\psi$ be a $\Delta$-linear support function with
$\psi(v)=\langle m_\sigma, v\rangle$ for each $v\subset\sigma\in\Delta$.
we will say that $\psi$ is {\em  $k$-convex} if for any
$\sigma\in\Delta$ and $v\not\subset\sigma$ $$\langle m_\sigma,v\rangle\geq \psi(v)+k
$$
\end{definition}
Confronting the notion of convex and strictly convex function, see \cite{Oda}, it is clear that
\begin{itemize}
\item $\psi$ is $0$-convex if and only if it is convex;
\item $\psi$ is $1$-convex if and only if it is strongly convex.
\end{itemize}
\rem\label{additive}It is clear from the definition that:
\begin{itemize}
\item If $\psi$ is $k$ convex then it is $t$-convex for any $t\leq k$;
\item If $\psi_1$ is $t_1$-convex and $\psi_2$ is $t_2$-convex, then $(\psi_1+\psi_2)$ is $(t_1+t_2)$-convex.
\end{itemize}
\vskip10pt
The meaning of convexity and strong convexity of a $\D$-support function associated to a line bundle $L$ is quite clear at least at the fixed points $x(\s)\in\D$. If $\psi_L$ is convex then:
$${\cx}^{m_{\s}}(x(\s))=\prod_1^n \cx_i^{a_1-a_i}\neq 0$$
If $\s=\la \rho_1,\cdots,\rho_n\ra$ and $\s'=\la \rho_0,\rho_2,\cdots,\rho_n\ra$ then
$${\cx}^{m_{\s'}}(x(\s))=\prod_2^n \cx_i^{a_1-a_i}\cx_1^{\la m_{\s'},\rho_1\ra+a_1}=0$$
in the case $\la m_{\s'},\rho_1\ra+a_1>0$, i.e. $\psi_L$ strictly convex.
In other words if $\psi_L$ is convex then for each invariant point there is a non vanishing section, and $\psi_L$ strictly convex implies that different invariant points can be separated.\\
The notion of $k$-convexity generalizes the above property to more points with possible multiplicities.

More geometrically a $\D$-support function $\psi$ is $k$-convex if
for each $\sigma\in\Delta$ the graph of the defining linear  function
$\la m_{\s},\,\,\ra$ is ``very" high compared to the graph of $\psi$.\\
Recall that if $\psi_L$ is convex then the polyhedra $P_L$ is the convex hull of the points $m_{\s}$ in $M_{\Bbb R}$. If $\psi_L$ is strictly convex then there is a correspondence between the faces in $\D$ and the set of non empty faces of $P_L$ (cf. \cite[2.12]{Oda}). Any face $F\subset P_L$ corresponds to
$$F^*=\{n\in N_{\Bbb R}|<m,n>=\psi_L(n),\text{ for any }m\in F\}\in\D$$
and any cone $\s\in\D$ corresponds to
$$\s^*=\{m\in M_{\Bbb R}|<m,n>=\psi_L(n),\text{ for any }n\in \s\}\subset P_L$$ Then $\psi_L$ being $k$-convex means that the length of the edges of the polyhedra, corresponding to $\tau=\s_i\cap\s_j$, are bigger or equal to $k$. This implies that $P_L$ is ``big enough" to choose points in it, corresponding to sections with an arbitrary prescribed jet.
\begin{lemma}\label{lenght}
Let $P_L$ be the polyhedra associated to $L$ and assume $\psi_L$ is $k$-convex, then
\begin{enumerate}
\item For each $\tau=\s_i\cap\s_j\in\D(n-1)$ the length of the associated edge $\tau^*$ is\\ $l(\tau^*)=|m_{\s_i}-m_{\s_j}|=l_{i,j}\geq k$;
\item Let $\s_1,\cdots,\s_r$ be the $n$-dimensional cones in $\D$. For each partition \\
$(t_1^1,\cdots,t_n^1,t_1^2,\cdots,t_n^2,\cdots,t_n^r)$ where $t_i^j\geq 0$, $\sum_{j=1}^nt_j^i=k_i-1$ and $\sum_1^rk_i=k+1$, and for any $\s_i\in\D(n)$ we can find $m\in P_L$ such that
\begin{itemize}
\item $\la m, \rho_l\ra=-a_l+t_l^i$ for all $\rho_j\subset\s_i$
\item $\la m, \rho_l\ra\geq-a_l+ t_l^j$ for all $\rho_j\not\subset\s_i$ with equality only if $t_l=0$
\end{itemize}
\end{enumerate}
\end{lemma}
\begin{pf} of (1). Let $\s_i=\la \rho_1,\cdots,\rho_n\ra$ and $\s_j=\la \rho_2,\cdots,\rho_{n+1}\ra$. Then using the basis $\{m_1,\cdots,m_n\}$ dual to $\{\rho_1,\cdots,\rho_n\}$
$m_{\s_i}=(-a_1,\cdots,-a_n)$ and $m_{\s_j}=(-a_1+l_{i,j},\cdots,-a_n)$, where $l_{i,j}\geq k$ since $\psi_L$ is $k$-convex.\end{pf}
\begin{pf} of (2). Fix $\s=\la \rho_1,\cdots,\rho_n\ra$ for simplicity of notation and let \\$\s_i=\la \rho_1,\cdots,\check{\rho_i},\cdots,\rho_n,\ol{\rho_i}\ra$ the $n$-cone so that $\s\cap\s_i=\tau_i=\la \rho_1,\cdots,\check{\rho_i},\cdots,\rho_n\ra$. By (1) the edge $\tau_i^*$ has lenght at least $k$. Choose the $t_i$-th lattice point next to $m_{\s}$ traveling on $\tau^*$ towards $m_{\s_i}$, i.e 
$$\ol{m_i}=(-a_1,\cdots,-a_i+t_i,\cdots,-a_n)=m_{\s}+(\frac {t_i}{l_i})(m_{\s_i}-m_{\s_i})$$ in the basis dual to $\{\rho_1,\cdots,\rho_n\}$, where $\la m_{\s_i},\rho_i\ra=-a_i+l_i\geq -a_i+k$ by hypothesis. Traveling on the $n$ edges next to $m_{\s}$ we get
$$m=m_{\s}+\sum_1^n(\frac {t_i}{l_i})(m_{\s_i}-m_{\s})$$
Rewriting it in the form
$$m=(1-\sum_1^n(\frac {t_i}{l_i}))m_{\s}+\sum_1^n(\frac {t_i}{l_i})m_{\s_i}$$
It is clear that $m$ is a convex combination of $\{m_{\s},m_{\s_1},\cdots,m_{\s_n}\}$, since \\$0\leq \sum_1^n(\frac {t_i}{l_i})\leq \frac{1}{k}\sum_1^n t_i\leq 1$ and therefore $m\in P_L=Conv(m_{\s})_{\s\in\D(n)}$. Moreover
\begin{itemize}
\item $\la m,\rho_i\ra=-a_i+t_i$ for $i=1,\cdots,n$
\item if $\rho_l\not\subset\s$ then
\begin{align*}
\la m,\rho_l\ra&=(1-\sum (\frac{t_i}{l_i}))\la m_{\s},\rho_l\ra+\sum [(\frac{t_i}{l_i})\la m_{\s_i},\rho_l\ra]\\
&> (1-\sum (\frac{t_i}{l_i}))(-a_l+k)+\sum [(\frac{t_i}{l_i})(-a_l)]\\
&=-a_l-k(\sum \frac{t_i}{l_i})+k\geq -a_l+(k-\sum t_i)
\end{align*}
If $k_j=0$ for all $j\neq i$ then $\sum t_i=k_i-1=k$ and thus $\la m,\rho_l\ra\geq -a_l.$ Otherwise $k-\sum t_i<t_l$ and thus $\la m,\rho_l\ra>-a_l+t_l.$
\end{itemize}
\end{pf}

\vskip10pt
\section{The Main Result}
Let us first formulate an equivalent criterion for $\psi_L$ to be $k$-convex in terms of the intersections of the divisor $L$ with the $T$-invariant rational curves associated to each $\tau\in\D(n-1)$.\\
 In fact the polyhedra $P_L$ having edges of length at least $k$ translates to the restriction of $L$ to each curve, corresponding to such edges, being at least $k$.\\ This is in a way a generalization of the ``toric Nakai criterion", cfr. \cite[Th. 2.18]{Oda}.
\begin{proposition}\label{inter} Let $L$ be a line bundle on a smooth $n$-dimensional toric variety $X$. Then $\psi_L$ is $k$-convex if and only if the restriction $L|_{V(\tau)}$ has degree $\geq k$, for every $\tau\in\D(n-1)$.\\
Equivalently if and only if $L_{V(\tau)}=\cO_{\pn{1}}(a)$ with $a\geq k$ for every $\tau\in\D(n-1)$. 
\end{proposition} 
\proof Let $\tau=\s_0\cap\s_1$ and assume $\s_i=\la \tau, n_i\ra$ for $i=0,1$. Then, since we are assuming $X$ to be non singular, there exists a $\Bbb Z$-basis
$\la n_i, n_2,\cdots,n_n\ra$ and $(n-1)$ integers $(s_2,\cdots,s_n)$ such that:
$$n_0+n_1-\sum_2^{n-2}s_in_i=0$$ 
Write $L=-\sum_i \psi_L(n_i)D_i$ where $D_1, D_0$ are the principal divisors associated to the edges $n_0, n_2$ and $D_i$ are the ones associated to $n_i$, $i=2,\cdots,n$.
then
\begin{itemize}
\item $D_1\cdot V(\tau)=D_0\cdot  V(\tau)=1$
\item $D_i\cdot  V(\tau)=-s_i$ for $i=2,\cdots,n$
\item $D_j\cdot  V(\tau)=0$ otherwise
\end{itemize}
\[\begin{array}{ll}
L\cdot V(\tau)&=\sum (-\psi_L(n_i))D_i\cdot V(\tau)=\\
&=-\psi_L(n_0)-\psi_L(n_1)+\sum_2^n \psi_L(n_i)s_i\\
&=\la m_{\s_1}, n_0\ra -\psi_L(n_0)
\end{array}\]
It follows that $L\cdot V(\tau)\geq k$ for all $\tau\in\D(n-1)$ if and only if for any $\s\in\D$ and $\la\rho_j, \s\cap\s'\ra=\s'$ the inequality 
$$\la m_{\s}, \rho_j\ra -\psi_L(\rho_j)\geq k$$
In other words $L\cdot V(\tau)\geq k$ if and only if the support function $\psi_L$ is $k$-convex.\qed

The relation between $k$-convex $\D$-support functions and $k$-very ampleness is given by the following:
\begin{theorem}\label{main} A line bundle  $L$ generates $k$-jets on $X$ if and only if the $\D$-support function $\psi_L$ is $k$-convex.
\end{theorem}
We need the following reduction step:\\
\noindent{\rm C\sc laim} {\em If $L$ is $k$-jet ample on any $r$-ple of fixed points $\{x(\s_1),\cdots,x(\s_r)\}$, $\s_i\neq \s_j$, then it is $k$-jet ample on $X$.}   
\begin{pf}{\rm \sc of the claim.} To every $(x_1,\cdots,x_{k+1})\in X^{k+1}$ we can associate $\cZ=(x_1,\cdots,x_r)$, a collection of $r\leq k+1$ distinct points of $X$, and $(k_1,\cdots,k_r)$, an $r$-ple of positive integers such that $\sum k_i=k+1$, simply counting the multeplicities of each $x_i$: 
$$(\underbrace{x_1,\cdots,x_1}_{k_1},\underbrace{x_2,\cdots,x_2}_{k_2},\cdots,\underbrace{x_r,\cdots,x_r}_{k_r})=[\cZ=(x_1,\cdots,x_r),(k_1,\cdots,k_r)]$$
Then to each $\underline{x}\in X^{k+1}$ we can associate the map:
$$\psi_{\underline{x}}=\psi_{\cZ}^{k_1,\cdots,k_r}:H^0(X,L)\to \bigoplus_1^r (J_{k_{i-1}}(L))_{x_i} $$
Let $C=\{\underline{x}\in X^{k+1}\text{ such that }coker(\psi_{\underline{x}})\neq 0\}$. Since $\psi_{\underline{x}}$ is an equivariant map, $C$ inherits the torus action from $X$, i.e. it is an invariant closed subvariety of $X^{k+1}$, and hence proper. If $L$ is not $k$-jet ample on $X$, then $\psi_{\underline{x}}$ is not surjective for some $\underline{x}=[\cZ,(k_1,\cdots,k_r)]$, which means $C\neq\emptyset$. But this implies $C^{T}\neq\emptyset$, where $C^{T}$ is the set of the fixed points in $C$. To see this one can apply Borel's fixed point theorem ( see \cite[21.1]{Hum}) or more directly observe that $C$ is a lower dimensional toric variety and thus it must contain fixed points. It follows that there exists $\underline{x}\in X^{k+1}$, fixed by the torus action, for which $\psi_{\underline{x}}$ is not surjective. Such $\underline{x}$ must have all the components fixed so it is of the form $(x(\s_1),\cdots,x(\s_r))$, which is a contradiction.\end{pf}
\begin{pf}{\sc of the theorem.}
If $L$ is a $k$-jet ample line bundle then the restriction $L|_{V(\tau)}$ to every $\tau\in\D(n-1)$ is $k$-jet ample, i.e. $L_{V(\tau)}=\cO_{\pn{1}}(a)$ with $a\geq k$ for every $\tau\in\D(n-1)$.
Proposition \ref{inter} then implies that $\psi_L$ is $k$-convex.\\ 
  Assume now that $\psi_L$ is $k$-convex. By the reduction step it suffices to prove that the map $\psi_{\cZ}^{(k_1,\cdots,k_r)}$ is surjective for each $\cZ=\{x(\s_1),\cdots,x(\s_r)\}$, with $k_1+\cdots+k_r=k+1$. This follows immediately from Lemma \ref{lenght}. For each $k_i$ and for each partition $t_1^i+\cdots+t_n^i=k_{i}-1$ we can choose $m\in P_L$ such that
\begin{itemize}
\item $\cx^{m}=\prod_{\rho_j\subset\s_i}\cx_j^{t_j^i}$ around $x(\s_i)$ and
\item $\cx^{m}=\prod_{\rho_j\subset\s_l}\cx_j^{t_j^l+c_j^l}$ around $x(\s_l)\neq x(\s_i)$, with $c_j^l>o$ for some $j$ and for any partition $\sum t_j^l\leq k_l$
\end{itemize}
This means that
$$(i_{k_i}( \cx^{m}(x(\s_1)),\cdots,i_{k_r}( \cx^{m}(x(\s_r)))=(0,\cdots,1,0,\cdots,0)$$
the non zero term corresponding to $\frac{\partial \cx^{m}}{\pr{x_1}{t_1}\cdots\pr{x_r}{t_r}}|_{x=x(\s_i)}$. Enough to prove the surjectivity.\end{pf}

\rem\label{add}From \ref{additive} it follows immediately that:
\begin{itemize}
\item if $L$ is $k$-jet ample then it is $t$-jet ample for any $t\leq k$;
\item If $L$ is $k$-jet ample  and $E$ is $t$-jet ample then the line bundle $E\otimes L$ is $(k+t)$-jet ample. In fact $\psi_{E\otimes L}=\psi_L+\psi_E$.
\end{itemize}
It is worth observing that in the toric case the notion of $k$-jet ampleness is equivalent to the notion of $k$-very ampleness (which is weaker in general).
For basic properties of $k$-very ample line bundles we refer to \cite{BeSoB}.
\begin{definition}
$L$ is said to be $k$-very ample if for every zeroscheme $(\cZ,\cO_{\cZ})$ of length $h^0(\cO_{\cZ})=k+1$ the map $H^0(X,L)\to H^0(\cZ,L\otimes\cO_{\cZ})$ is onto.
\end{definition}
\begin{proposition}
A line bundle $L$ on a  smooth toric variety $X$ is $k$-very ample if and only if it is $k$-jet ample.
\end{proposition}
\proof If $L$ is $k$-jet ample then it is $k$-very ample ( see \cite[Prop. 2.2]{BeSok}).\\
If $L$ is $k$-very ample then the degree of $L$ restricted to any irreducible curve must be $\geq k$, i.e. $L\cdot V(\tau)\geq k$ for all $\tau\in\D(n-1)$ and thus it is $k$-jet ample by \ref{inter}. \qed
\section{Examples}
In this section we work out few examples, for which the $k$-jet ampleness has been studied otherwise, to convince the reader that this is in fact the right way to formulate the result.
\ex\label{pn}{\sc The projective space $\pn{n}$.} Notation as in \ref{pn1}.
Let  $L=t_1D_1+\cdots+t_{n+1}D_{n+1}\cong (t_1+\cdots+t_{n+1})D_1$. By \ref{inter} $L$ is $k$-jet ample if and only if $V({\s_i})\cdot L\geq k$ where $\s_i$ is the n-dimensional cone $\la \rho_1,\cdots\check{\rho_i},\cdots,\rho_{n+1}\ra$ and $D_i$ is the divisor associated to the edge $\rho_i$. Since $\rho_{n+1}+\rho_i+\sum_{j\neq i}\rho_j=0$
$$V({\s_i})\cdot L=t_1+\cdot+t_{n+1}\geq k$$
In other words $\cO_{\pn{N}}(a)$ is $k$-jet ample if and only if it is $k$-very ample if and only if $a\geq k$, as proven in \cite{BeSok}.
\vskip10pt
\ex\label{fn}{\sc The Hirzebruch surface ${\Bbb F}_n$.}
Let $\{e_1,e_2\}$ be  the standard basis for ${\Bbb R}^2$. The Hirzebruch surface ${\Bbb F}_n$ is the toric surface associated to the fan $\D$ spanned by the following $2$-cones:
$$\s_1=\la e_1,e_2\ra,\,\s_2=\la e_2,-e_1\ra,\,\s_3=\la -e_2,-e_1+ne_2\ra,\,\s_4=\la -e_1+ne_2,e_2\ra$$
Let $D_1,\cdots,D_4$ be the divisors associated respectively to $e_2,e_1,-e_2,-e_1+ne_2$. Then
we have the following intersection matrix:
\[\begin{pmatrix}D_1^2&D_1\cdot D_2&D_1\cdot D_3&D_1\cdot D_4\cr
D_2\cdot D_1&D_2^2&D_2\cdot D_3&D_2\cdot D_4\cr
D_3\cdot D_1&D_2\cdot D_3&D_3^2&D_3\cdot D_4\cr
D_4\cdot D_1&D_2\cdot D_4&D_4\cdot D_3&D_4^2\cr
\end{pmatrix}=
\begin{pmatrix}-n&1&0&1\cr
         1&0&1&0\cr
         0&1&n&1\cr
         1&0&1&0\cr
\end{pmatrix}\]
Recall that $D_3\equiv D_1+nD_2$ and $D_2\equiv D_4$.
Let $L=a_1D_1+\cdots+a_2D_4=(a_1+a_3)D_1+(a_4+a_2+na_3)D_2$. Then  by \ref{inter} $L$ is $k$-jet ample if and only if it is $k$-very ample if and only if
\begin{itemize}
\item $L\cdot D_1=a_4-na_1+a_2\geq k$
\item $L\cdot D_2=a_1+a_3\geq k$
\item $L\cdot D_3=a_2+na_3+a_4\geq k$
\item $L\cdot D_4=a_1+a_2\geq k$
\end{itemize}
which of course is equivalent to saying that $L=aE_o+bf$, where $E_0$ is the section of minimal selfintersection $-n$ ( i.e. $D_1$) and $f$ the general fiber of the projection onto $\pn{1}$ (i.e. $D_2$), is $k$-jet ample if and only if
 \begin{itemize}
\item $a=a_1+a_3\geq k$ 
\item $-an+b=-na_1-na_3+a_4+na_3+a_2=a_4-na_1+a_2\geq k$
\end{itemize}
This conditions have been given by Beltrametti-Sommese in \cite{BeSok}, using a decomposition argument.
\vskip10pt
\ex{\sc Del Pezzo surfaces.}
The toric Del Pezzo surfaces are $\pn{2}$, ${\Bbb F}_1$ and the equivariant blow up of $\pn{2}$ in $2$ or $3$ points.
The most interesting one is the last one.
Let $S$ be the equivariant blow up of $\pn{2}$ in the $3$ invariant points as described in \ref{delpezzo}.
The principal divisors $D_1,\cdots,D_6$ are the $6$ $(-1)$-curves on the surfaces, i.e. the three exceptional divisors and the pull back of the three lines passing through two of the $3$ points blown up.\\
Proposition \ref{inter} says that $L$ is $k$-jet ample if and only if it is $k$-very ample if and only if the intersection with all the $(-1)$-curves on the surface is $\geq k$.\\
This criterion has been given for $k$-very ampleness in \cite{DR} using a generalization of Reider's theorem.\\
Note that the equivalence between $k$-jet ampleness and $k$-very ampleness is not always true for Del Pezzo surfaces. In fact it is not hard to see that if $S$ is the blow up of $\pn{2}$ in $7$ points in general position ( so it is not toric), then the line bundle $L=-2K_S$ is $2$-very ample but it is not $2$-jet ample.
\section{ Local positivity applications}\label{applications}
In this section we report some nice applications of $k$-jet ampleness to the study of ``local positivity" of line bundles. Most of it is a survey on well known results. We think it is interesting to show how these results can be established rather easily in the case of toric varieties.
\vskip10pt
\noindent{\sc Blow ups}. A $k$-jet ample line bundle carries its positivity along blow ups at a finite number of points. The following property has been proved by Beltrametti and Sommese in \cite{BeSo96} and it has a very ``visible" proof in the toric case. We refer to \cite{Oda} for notation and definition of equivariant blow ups.
\begin{proposition}\label{blowup} Let $p:X(\D')\to X(\D)$ be the equivariant blow up of $X(\D)$ at $r$ points, $x_1,\cdots,x_r$, and let $L$ be a $k$-jet ample line bundle on $X(\D)$. Then $p^*(L)-\sum\e_iE_i$ is min$(k-\sum\e_i,\e_1,\cdots,\e_r)$-jet ample on $X(\D)$, where $E_i$'s are the exceptional divisors.
\end{proposition}
\begin{pf}Use induction on $r$. Assume the number of edges in $\D$ is $N$. \\
If $r=1$, let $x=x(\s)$, $\s=\la \rho_1,\cdots,\rho_n\ra$ and $E_1=E$ be the divisor associated to the edge $\overline{\rho}=\rho_1+\cdots+\rho_n$. In the new fan $\D'$ there are $n$ new $n$-cones $\s_i=\la \overline{\rho},\rho_1,\cdots,\check{\rho_i},\cdots,\rho_n\ra$. Moreover let $\overline{D_i}$ be the divisors in $Pic(X(\D'))$ corresponding to the edges $\rho_i$, then 
\begin{itemize}
\item $\overline{D_i}=p^*(D_i)-E$ for $i=1,...,n$;
\item $\overline{D_i}=D_i$ for $i=n+1,...,N$.
\end{itemize}
 If $L=\sum a_i D_i$ then
$$\cH=p^*(L)-\e E=\sum_1^N a_ip^*(D_i)-\e E=\sum_1^N a_i\overline{D_i}-(\e+\sum_i^n a_i)E$$
Let $\tau\in\D'(n-1)$. If $\tau\in\D(n-1)$ then clearly $\cH\cdot V(\tau)=L\cdot V(\tau)\geq k$. If $\tau\in\D'-\D$ then it is one of the following:
\begin{itemize}
\item[(a)] $\s_i\cap\s_j=\la \overline{\rho},\rho_1,\cdots,\check{\rho_i},\check{\rho_j},\cdots,\rho_n\ra$
\item[(b)] $\s_i\cap\la\rho_{n+1},\rho_1,\cdots,\check{\rho_i},\cdots,\rho_n\ra$
\end{itemize}
Following the lines of \ref{inter}:\\
In case (a), since $\rho_i+\rho_j-\overline{\rho}+\sum_{l\neq i,j}\rho_l=0$
$$\cH\cdot V(\tau)=-a_1-a_2+\e+\sum_1^n a_l+
\sum_{l\neq i,j}a_l=\e$$
In case (b), assume $\rho_{n+1}+\rho_i-\sum_{j\neq i}s_j\rho_j=0$, then $\rho_n+1+\overline{\rho}-\sum_{j\neq i}(s_j-1)\rho_j=0$ and
$$\cH\cdot V(\tau)=-a_{n+1}-\e-\sum_1^n a_j+\sum_{j\neq i}s_j+\sum_{j\neq i}a_i=L\cdot V(\tau')-\e$$
where $\tau'=\la \rho_1,\cdots,\check{\rho_i},\cdots,\rho_n,\rho_{n+1}\ra\cap\s\in\D(n-1)$\\
If $r>1$, iterating this process, after $r$ blow-ups $\cH=p_{r}^*(L')-\e_rE$, where \\$p^*_r:X(\D')\to X(\D_{r-1})$ is the $r$-th blow up map. By induction $L'$ is min$(k-\sum_1^{r-1}\e_1,\e_r,\cdots,\e_{r-1})$-jet ample on $X(\D_{r-1})$. Clearely from what done before 
$$\cH\cdot V(\tau)\geq min(L'\cdot V(\tau')-\e_r,\e_r)\geq min(k-\sum_1^n\e_i,\e_1,\cdots,\e_r)$$ for any $\tau\in \D'(n-1)$.\end{pf}
\vskip10pt
\noindent{\sc Toric Seshadri criterion}. An ample line bundle on a smooth projective variety, $X$, is characterized by the positive value of its Seshadri constant at each point.\\
Let $L$ be a nef line bundle on $X$. For every irreducible curve $C\subset X$,  $m_x(C)$ denotes the multiplicity of $C$ at the point $x\in C$ and
$$m(C)=\sup_{x\in C}\{m_x(C)\}$$.
\begin{theorem}(Seshadri \cite[7.1]{Ha}) A line bundle $L$ on $X$ is ample if and only if there exists $\e>0$ such that $L\cdot C\geq \e\cdot m(C)$ for every irreducible curve $C\subset X$.
\end{theorem}
As for the Nakai criterion we can generalize the Seshadri criterion to $k$-jet ampleness on toric varieties.\\
Let us first reformulate the Seshadri's theorem in the ``modern languige" of Seshadri constants.\\
For a nef line bundle $L$ the Seshadri constant of $L$ at a point $x\in X$ is the real number
$$\e(L,x)=\inf_{x\in C}\frac{L\cdot C}{m_x(C)}=\sup_{\e}\{\e\in {\Bbb R}|p^*(L)-\e L\text{ is nef }\}$$
where the inf is taken over all the irreducible curves containing $x$ and $p$ is the blow-up map of $X$ at $x$.
Then one can immediately see that The Seshadri criterion says that $L$ is ample if and only if $\e(L,x)>0$ for every $x\in X$.\\
Demailly showed that the Seshadri constant is a measure of the highest degree jets that can be generated by the global sections of $L$:
\begin{proposition}\label{dem}\cite{Dem}
Let $s(L,x)$ be the largest integer such that $|L|$ generates $s$-jets at $x$. Then
$$\e(L,x)= \limsup_{n\to\infty}\frac{s(nL,x)}{n}$$
\end{proposition}
Since every ample line bundle on a toric variety is very ample the toric Seshadri criterion says that $L$ is a ample if and only if $\e(L,x)\geq 1$ for every $x\in X(\D)$. More generally:
\begin{proposition}\label{ses} 
A line bundle $L$ on a non singular toric variety $X(\D)$ is $k$-jet ample if and only if there exists an $\e\geq k$ such that 
$$L\cdot V(\tau)\geq \e \cdot m(V(\tau))$$ for every invariant curve $V(\tau)$.\end{proposition}
\begin{pf} Assume $L$ is $k$-jet ample, where $k$ is the biggest integer such that the property is true. Then by \ref{dem} $\e(L,x)=k$, since $nL$ is $(nk)$-jet ample by \ref{add}. It follows that $\frac{L\cdot V(\tau)}{m(V(\tau))}\geq k$ for every $\tau\in\D$.\\
Assume now that $L\cdot V(\tau)\geq k\cdot m(V(\tau))$ for any $\tau\in\D$. Let $\tau=\s_i\cap\s_j$, we have to prove that $L\cdot V(\tau)\geq k$.
Consider the invariant point $x=V(\s_i)\in V(\tau)$, clearely $m_{x}(V(\tau))=1$. Then
$L\cdot V(\tau)\geq k\cdot m(V(\tau))\geq k\cdot m_x(V(\tau))\geq k$.
\end{pf}
\begin{corollary}\label{sesh} A line bundle $L$ on a non singular toric variety is $k$-jet ample if and only if $\e(L,x)\geq k$ for every $x\in X$.\end{corollary}
\begin{pf} If $L$ is $k$-jet ample then $\e(L,x)=k$ by \ref{dem}.\\
If $\e(L,x)\geq k$ for every $x\in X$ then in particular $L\cdot V(\tau)\geq k \cdot m(V(\tau))$ for every invariant curve $V(\tau)$ and thus $L$ is $k$-jet ample by \ref{ses}.\end{pf}
\vskip10pt
\noindent{\sc Higher adjoint series}. By use of bounds on the Seshadry constant of $L$ at a sufficiently general point Ein and Lazarsfeld showed that:
\begin{proposition}\cite[5.14]{Laz} Let $L$ be an ample line bundle on a smooth surface $S$. Then the adjoint series $|K_S+(k+3)L|$ generates $k$-jets at a sufficiently general point $x\in S$.\end{proposition}
The fact that for a line bundle on a toric variety being ample is equivalent to being very ample implies a simple generalization.
\begin{proposition} Let $S$ be a non singular toric surface and let $L$ an ample line bundle on it such that $L^2>1$. Then
\begin{itemize}
\item[(a)] $|K_S+(k+2)L|$ generates $k$-jets at every point $x\in S$;
\item[(b)] $|K_S+(2k+2)|$ generates $k$-jets on $S$.
\end{itemize}
\end{proposition}
\begin{pf}(a) By the Nakai toric criterium $L$ is in fact very ample and by \ref{add} $(k+2)L$ is $(k+2)$-jet ample. Let $x\in S$, using the long exact sequence:
$$\to H^0(K_S+(k+2)L)\to H^0((K_S+(k+2)L)/{\frak m}_x^{k+1}))\to H^1((K_S+(k+2)L)\otimes {\frak m}_x^{k+1})\to$$
the vanishing of  $H^1((K_S+(k+2)L)\otimes {\frak m}_x^{k+1})$ would imply the result. Let $p:\ol{S}\to S$ be the blow up of $S$ at $x$ with $E=p^{-1}(x)$, then by Leray spectral sequence and Serre duality 
$$H^1((K_S+(k+2)L)\otimes {\frak m}_x^{k+1})=H^1(K_{\ol{S}}+[p^*((k+2)L)-(k+2)E])$$
Kawamata vanishing theorem applies since $p^*((k+2)L)-(k+2)E$ is nef and big and thus $K_S+(k+2)L$ is $k$-jet ample at any point $x\in S$.\\
(b) Consider now simoultaneus jets supported on $\{x_1,\cdots x_r\}\in S$ and $(k_1,\cdots,k_r)\in{\Bbb Z}_+^r$ such that $\sum k_1=k+1$. Let $p:\ol{S}\to S$ the blow up of $S$ at $x_1,\cdots, x_r$ with $E_i=p^{-1}(x_i)$.  Using the same exact sequence as above it suffices to prove that:
$$H^1((K_S+(2k+2)L)\otimes({\frak m}_{x_1}^{k_1}\otimes\cdots\otimes {\frak m}_{x_r}^{k_r}))=H^1(K_{\ol{S}}+[p^*((2k+2)L)-\sum (k_i+1)E_i])=0$$
Since $(2k+2)L$ is $(2k+2)$-jet ample and $\sum(k_i+1)\leq 2k+2$, $p^*((2k+2)L)-\sum (k_i+1)E_i$ is spanned by \ref{blowup}. Moreover
$$(p^*((2k+2)L)-\sum (k_i+1)E_i)^2>(2k+2-\sum(k_1+1))(2k+2+\sum(k_1+1))\geq 0$$
Then Kawamata vanishing theorem applies to gine the nedeed vanishing.\end{pf}
\vskip10pt
\noindent{\sc The $k$-reduction}. In the case of surfaces we can make some further remarks about the ``$k$-reduction" process (see \cite{BeSo}).
Let $S$ be a non singular toric surface and $L$ a $k$-very ample line bundle on it. Since $k$-jet very ampleness and $k$-very ampleness are equivalent we will use freely the property of being $k$-very ample according to the criterion given in \ref{main}.
Assume the adjoint bundle is not $k$-very ample i.e., the surface contains $(-1)$-curves whose intersection with
$L$ is exactly $k$. If the $k$-adjoint bundle
$kK_S+ L$ is nef we can contract down those curves and get the 
$k$-reduction $(S',L')$.

Notice that if $L=-kK_{S}$ then $-K_{S}$ is ample and hence $S$ is a toric
Del Pezzo surface. We can then compute directly the nefness of the $k$-adjoint bundle obtaining the same result as in \cite{BeSo}. Notation as in \ref{pn},\ref{fn}. 
\begin{proposition} Let $L\neq -kK_{S}$ be a $k$-very ample line bundle on
$S$. Then $kK_{S}+L$ is nef unless:
\begin{itemize}
\item $S=\pn{2}$ and $L=aD_1$ with $a<3k$ ; 
\item $S={\Bbb F}_r$ and $L=aD_1+bD_2$ with $a< 2k$.
\end{itemize}
\end{proposition}
\begin{pf}Let $L=\sum a_iD_i$ and $D_i^2=-s_i$, then
$$(kK_S+L)\cdot D_i=L\cdot D_i+k(s_i-2)\geq 0\text{ if }s_i\geq 1$$
 Recall that $S$ is isomorphic to $\pn{2}$, ${\Bbb F}_n$ or their equivariant blow up in a finite number of points. If $S$ is minimal then intersecting $L$ with the basic generators of $Pic(S)$ and imposing at least one intersection to be less then $k$ gives the cases in the statement. Assume now $S$ not minimal.
If $S=Bl_r(\pn{2})$ (i.e. the blow up of $\pn{2}$ in $r$ points) let $D_1,D_j,D_l$ be the divisors associated to the edges $(0,1),(1,0),(-1,-1)$ respectively. If $r\geq 2$, for each $D_i$ generator of $Pic(S)$, the corresponding weight $D_i^2=-s_i\leq -1$ unless possibly only one among $(s_1,s_j,s_l)$, say $s_1=-1$. But in this case $D_1\equiv D_2+\sum_1^rD_i$, the $D_i$'s being the exceptional divisors, $L\cdot D_1\geq (r+1)k$ and
$$(kK_S+L)D_1\geq (r+1)k-3k\geq 0$$
If $r=1$ then $S\cong {\Bbb F}_1$.\\
If $S=Bl_r({\Bbb F}_n)$, let $D_l,D_j,D_3,D_h$ be the divisors corresponding to the edges $(0,1),\\
(1,0),(0,-1),(-1,n)$ respectively. We can assume the weights $s_i\geq 1 $, unless possibly $s_3=-n+s<1$, since $Bl_1{\Bbb F}_0\cong Bl_2(\pn{2})$ and $D_j\equiv -D_h$. But in this case $D_3\equiv D_l+nD_j+\sum_1^{r-s}D_i$, $L\cdot D_3\geq (n+1+r-s)k$ and
$$(kK_S+L)D_3\geq (n+r+1-s)k-(n-s+2)k=(r-1)k\geq 0$$\end{pf}


\end{document}